\documentclass[10pt,twoside]{article}
\usepackage{amssymb}
\usepackage{latexsym}
\usepackage{fullpage}
\usepackage{latexsym,rotating,amsmath,epsfig}
\usepackage{here}
\usepackage{amsmath}
\usepackage{amscd}
\usepackage{oldgerm}
\usepackage{fancyheadings}

\def\N{\mathbb{N}}

\def\R{\mathbb{R}}

\def\1{\mathbf{1}}
\def\:{\lrcorner}
\def\#{\sharp}

\def\l{\lambda}
\def\a{\alpha}

\def\g{\gamma}
\def\d{\delta}
\def\e{\epsilon}

\def\o{\circ}

\def\<#1,#2>{\langle#1,\,#2\rangle}

\def\K{\mathbb{K}\,}

\def\qed{\ensuremath{\quad\Box\quad}}
\def\pfill{\par\vskip2mm plus1mm minus1mm\noindent}
\def\inv#1{\raise.1em\hbox to 0pt{$^{-1}$\hss}_{#1}\;}

\def\v{\noindent}

\def\ov{\overline}
\newcommand{\bean}{\begin{eqnarray*}}
\newcommand{\eean}{\end{eqnarray*}}
\newcommand{\benu}{\begin{enumerate}}
\newcommand{\eenu}{\end{enumerate}}
\newcommand{\eea}{\end{eqnarray}}
\newcommand{\bea}{\begin{eqnarray}}

\newtheorem{Theorem}{Theorem}[section]
\newtheorem{Lemma}[Theorem]{Lemma}

\newtheorem{Definition}[Theorem]{Definition}

\newtheorem{Proposition}[Theorem]{Proposition}

\title{Metrizability of spaces of homomorphisms between metric vector spaces}
\begin{document}

\author{Olaf M\"uller\footnote{Instituto de Matem\'aticas, Universidad Nacional Aut\'onoma de M\'exico (UNAM) Campus Morelia, C. P. 58190, Morelia, Michoac\'an, Mexico. email: olaf@matmor.unam.mx}}

\maketitle

\begin{abstract}
\v This note tries to give an answer to the following question: Is there a sufficiently rich class of metric vector spaces such that sufficiently large spaces of continuous linear maps between them are metrizable? 
\end{abstract}

\section{Introduction}

\v A useful feature of finite-dimensional analysis is the representation of the differential of a $C^1$ map $f: \R^n \supset U \rightarrow \R^m$ as a map taking values in the metric space $CL(\R^n, \R^m)$. In Banach spaces this point of view begins to cause problems and leads to the distinction between the (weaker) notion of Michal-Bastiani differentiability and the (stronger) notion of Fr\'echet differentiability (cf. \cite{hG}, where it is shown that the two notions differ only by one degree of differentiability). As soon as leaving normable spaces towards more general metric vector spaces, the situation gets even more complicated. On the other hand, it is well-known that we cannot avoid non-normable spaces if we want to include spaces in which derivative operators are continuous: 

\begin{Theorem}
There is no normable topology on $F_0 := C^{\infty} ([0,1], \R)$ with continuous derivative operator.
\end{Theorem}

\v {\bf Proof.} This is because continuity in normed spaces is equivalent to boundedness but there are arbitrarily high eigenvalues of the differential operator $\partial_t$ given by the functions $t \mapsto sin (Nt)$. \hfill \qed

\bigskip




\v Now, most people working in Fr\'echet spaces tend to consider only Michal-Bastiani differentiability which renouns completely on the concept of the differential taking values in spaces of linear maps and considering it as a map $f': U \times V \rightarrow W$. One reason for this is that if one equips the space of continuous linear maps $CL(F,G)$ between two Fr\'echet spaces with a topology, it turns out that in very general cases the evaluation map is not continuous any more (for a good overview cf. \cite{hhK} who suggested to circumvent this problem by considering convergence structures instead of topologies). But in the light of recent results on inverse function theorems for so-called bounded differentiable maps (\cite{oM}) it seems desirable to explore other types of differentiability which do include some form of iterated spaces of homomorphisms of the type above. The program followed by this and subsequent notes will thus be to find the exact reasons of non-metrizability and to provide some appropriate classes of metric vector spaces and linear maps between them.

\v The research leading to this note has been partially funded by the CONACyT project 82471.

\section{Palette topologies on $CL(V,W)$}

\v Often topologies on the dual space of a tvs $V$ are defined by means of a family of subspaces of $V$. Here and subsequently $CL(V,W)$ is the space of continuous linear maps from a tvs $V$ to a tvs $W$, and for $A \subset V$, $B \subset W$, $(A,B)$ is defined as the subset of $CL(V,W)$ which consists of all maps that map $A$ into $B$. For a family $F$ of subsets of $V$ the topology $\tau_F$ on $CL(V,W)$ is defined as the topology generated by the sets $(P, O)$ as a subbasis where $P \in F$ and $O \subset W$ open. It turns out that some properties of $F$ are important to ensure that $\tau_F$ is a vector space topology.

\begin{Definition}
Let $V,W$ be lhs, $A$ a linear subspace of the continuous linear maps from $V$ to $W$. An {\bf $A$-palette} is a subset $P$ of the potence set $\mathcal{P}(V)$ of $V$ with the properties

\begin{enumerate}
\item{For every member $p$ of $P$ and every member $a$ of $A$, $a(p)$ is s-bounded in $W$,}
\item{$A, B  \in P \Rightarrow A \cup B \in P$,}
\item{$A, B \in P \Rightarrow \l A  \in P$,}
\item{$A \in P, v \in V  \Rightarrow \ov{conv} (A, \{ v \}) \in P$,}
\item{$\bigcup_{A \in P} A$ is dense in $  V$.} 
\end{enumerate}

\v The palette is called {\bf strong} iff, for every neighborhood $N$ of $0$ it contains an element $P_N$ with $P_N \subset N$. 
\end{Definition}

\begin{Theorem}[for (i) cf. \cite{hS} III.3.1 and III.3.2]
Let $V$ be an lhs. 

\v (i) Every $A$-palette $P$ of $V$ generates as a subbasis a lhs topology $\tau_P$ on $A \subset CL(V,W)$. 

\v (ii) The evaluation map $ev_x$ for $x \in V$ is continuous on $(A, \tau_P)$ if and only if $P$ is strong.
\end{Theorem}

\v {\bf Proof.} (i) We have to show that there is a basis $B$ of zero neighborhoods which are point-absorbing, circled and such that for every $W \in B$ there is an $U \in B$ with $U+U \subset W$. If we choose a zero neighborhood base $H$ of $W$ which consists of circled sets, then the sets $(S,h)$ with $h \in H$ are circled as well as $\l (S,h) = (S, \l h)$ (and even convex, so the topology defined is locally convex). Also the last property is satisfied automatically. The first property is equivalent to the first property in the definition of palettes. It remains to show that $\tau_P$ is Hausdorff. Thus we have to show that $0 \in CL(V,W)$ is closed. Let $0 \neq f \in CL(V,W)$, then there is a $x \in V$ with $f(x) \neq 0$. Hausdorffness of $W$ allows us to find a zero neighborhood $W$ with $0 \notin f(x) + W$ and some zero neighborhood $U$ with $U + U \subset W$. By continuity of $f$ and density of $P$ we find a $p_0 \in P$ with $f(p_0)  \cap (f(x) + U) \neq \emptyset $, and then $0 \notin (p_0 , f(x) + U) \ni f  $.

\v (ii) Now we want to show that if the palette is strong, the evaluation map is continuous: We want to show openness of $ev_x^{-1} (U) = (x, U)$ for an open set $U$ in $W$ and every $x \in V$. If $L \in (x, U)$ then, as $\cdot - L$ is continuous, the subset $(x,U)$ is open iff $(x, U)-L= (x, J:= U- L(x))$ is open, and $J$ is an open neighborhood of $0$. But because $A$ consists of continuous maps, $(x,J) = \bigcup_{N \in N(0)}(x+N ,J)$. Now because of the strongness condition $ P_N \subset N $ we get $(x,J) = \bigcup_{p \in P^0} (x+ p, J) $ where $P^0$ is the subfamily of $P$ consisting of the elements containing $0$. But the sets $x + C_p p$ are members of the palette because of the defining properties 3 and 4, so $(x,J)$ is open and therefore $(x,U)$ as well. \hfill \qed     

\bigskip

\v {\bf Examples:} Basic examples are
\begin{enumerate}
\item{the palette $\mathcal{FC}$ of convex compact subsets contained in finite-dimensional linear subspaces,} 
\item{the palette $\mathcal{F}$ of compact subsets contained in finite-dimensional linear subspaces,} 
\item{the palette $\mathcal{CC}$ of convex compact subspaces} 
\item{the palette $\mathcal{C}$ of compact subspaces} 
\item{the palette $\mathcal{PC}$ of precompact subsets}
\item{the palette $\mathcal{S}$ of s-bounded subsets}
\item{the palette $\mathcal{B}_s $ of metrically bounded subsets of diameter $\leq s$}
\item{the palette $\mathcal{B} $ of metrically bounded subsets}
\end{enumerate}

\bigskip

\v {\bf Remark:} It is obvious that if $P_1 \subset P_2$ then $\tau_{P_1} \subset \tau_{P_2}$. 

\v {\bf Remark:} If the metric is bounded, then $\mathcal{B}$ is the maximal palette consisting of the whole potence set of $V$, and then it is easy to see that this generates a completely disconnected topology. 


\v {\bf Remark:} Property 5 is needed only to show Hausdorffness. Thus, for example, $B_rL(V,W)$ as described in \cite{oM} does not come from any palette, even without 5, as it is Hausdorff anyway.


\begin{Definition}
A {\bf fundamental sequence} of a palette $P$ is an increasing sequence $S_1 \subset S_2 \subset ...$ of elements of $P$ such that every element of $P$ is contained in some $P_k$.
\end{Definition}

\begin{Lemma}[\cite{gK}]
If $V$ is a metrizable tvs with a palette $P$ and an associated fundamental sequence $\{ S_i \vert i \in \N \}  $, then there is an $i \in \N$ such that $S_i$ absorbs every element of $P$. 
\end{Lemma}

\v {\bf Proof.} Suppose there is no such $S_i$. Then w.r.o.g. let $S_i$ be absolutely convex and such that no $S_n $ absorbs $S_{n+1}$. Choose a sequence $x_n \in S_1 \setminus \{ 0 \}$ with $x_n \rightarrow 0 $. Then for all $(n,k) \in \N^2 $ choose $z_{n,k} \in k^{-1} S_n$ with 

\bea
\label{nichtdrin}
z_{n,k} \notin (k+1) S_{n-1} .
\eea

\v Now define $M:= \{ x_n + z_{n,k} \vert n,k \in \N \}$. We want to show that the sequential completion $\ov{M}^s$ of $M$ is not closed in $V$, in contradiction to the assumption of metricity of $V$ which would define a metric on $M$ by restriction. To see that $\ov{M}^s$ is not closed, observe that $x_n \in \ov{M}^s$ and $x_n \rightarrow 0 $, but $0 \notin \ov{M}^s$: Suppose a sequence in $M$ converges to $0$, then it is s-bounded, so by the defining property of fundamental sequences it is contained in some $S_m$. Then we have $x_n + z_{n,k} \in S_m$ and therefore $z_{n,k} \in S_m - S_1 = S_m + S_1 \subset 2 S_m$, therefore because of Eq. \ref{nichtdrin} we have $n \leq m$ for all members of the sequence. But as the sequence converges to $0$, it has to contain arbitrarily high values of $n$, a contradiction.  \hfill \qed 

\begin{Theorem}[cf. \cite{wR}]
Every metrizable tvs has a translation-invariant compatible metric with circled balls. \hfill \qed
\end{Theorem}

\begin{Theorem}
\label{scalbound}
Let $(V,d)$ be a metric vector space with circled balls. Then $(V,d)$ is scalar-bounded by $2$, and we have even $d(sv, 0) \leq (-[-s]) d(v,0)$ for any $s>0$ where $[\cdot]$ is the Gauss bracket.
\end{Theorem}

\v {\bf Proof.} This is just an easy application of the triangle inequality. \hfill \qed 

\begin{Theorem}
\label{shithappens}
Let $V$ be a metrizable locally convex tvs, $A \subset V^* := CL(V, \R)$ with the Hahn-Banach property and an $A$-palette $P$ which contains the palette of convex compact sets. If $(A, \tau_P)$ is metrizable then there is an element of $P$ which contains an open set. 
\end{Theorem}

\v {\bf Remark.} Here we could also consider the space of all continuous linear maps into a Banach space. The latter property we will need in the rescaling process below.

\v {\bf Proof.} If $(A, \tau_P)$ is metrizable, then there is a countable system of zero neighborhoods $U_i$ which are w.r.o.g. of the form $U_i := (S_i, O_i)$ where $S_i \in P $ and $O_i \in \K$ (here we need the first property in the definition of palettes, stability under finite union, and because of $(S_1, O_1) \cap ... \cap (S_i , O_i) \supset (S_1 \cup ... S_i , O_1 \cap ... \cap O_i)$). By rescaling we can even find a system of the form $U_i' := (S_i' := k_i \cdot S_i, O_1 )$ as there are constants $k_i >0$ with $O_i \subset k_i O_1 $ by local s-boundedness of $\K$ and because of the second property in the definition of palettes. Then the $S_i'$ are a fundamental system for $P$, because otherwise there is $S \in P$ with $S \not\subset S_i'$ and therefore there is an open set $(S,O_1)$ not containing any $(S_i', O_1)$ (applying the Hahn-Banach property assumed above). With the lemma above we conclude that there is an $S_m'$ which absorbs all sets of $P$, thus all convex compact sets. Then $S_m'$ contains a ball: If not, define a sequence $x_n \in V \setminus S_m'$ but with $a_n := d(x_n, 0) \rightarrow 0$. Then $x_n' := \sqrt{a_n}^{-1} x_n \rightarrow 0$ because $(V,d)$ is scalar-bounded by $2$, therefore $conv ( \{ x_n \vert n \in \N \} )$ is convex and compact, but is not absorbed by $S_m'$ (as the necessary scaling factor to absorb the $n$-th point of the sequence would have to be smaller than $\sqrt{a_n}$ which tends to $0$), a contradiction. \hfill \qed

\bigskip

\v As a corollary we get the well-known result

\begin{Theorem} 
\label{et-tu-dual?}
Let $V$ be a metrizable locally convex tvs. If $V'$ equipped with the compact-open topology is metrizable, $V$ is finite-dimensional. If the s-bounded-open topology is metrizable, then $V$ is normable. \hfill \qed
\end{Theorem}

\v So let's go on with our quest: First let us look below the CCO topology. What about the FO (finite-open) topology on $CL(V,W)$ described by the palettes $\mathcal{FC}$ or $\mathcal{F}$? This topology is described by uniform convergence of filters on finite sets (as the maps are linear, finite and finite-dimensional are equivalent here) and is therefore complete if and only if $W$ is complete.

\begin{Theorem}
If $V$ is infinite-dimensional, the finite-open topology on $V^*$ is not metrizable.
\end{Theorem}

\v {\bf Proof.} If we assume that there is a countable zero neighborhood base $U1 \supset U_2 ...$ we can assume w.r.o.g. that $U_i = (F_i, O_i)$ with $F_i $ finite sets. This gives us a countable Hamel generating system (and by the usual clean-up procedure a countable Hamel basis) which does not exist in infinite-dimensional complete metric vector spaces as they are nonmeager in itself and as the sequence of finite-dimensional subspaces following the basis would be sequence of closed subspaces of empty interior whose union is the whole space. \hfill \qed

\section{General obstructions against metrization}

\v Now we will see that typical Fr\'echet spaces do not have well-behaved metrizable topologies on their dual spaces. Throughout this section, let $V_1:= \R^{\N}$ be the space of real sequences equipped with any vector space topology $\tau_1$ in which the linear maps $d_n \in L(V_1, \R) $ given by $d_n (a) := a_n$ are continuous. Let $V_2:= C^{\infty} ([0,1], \R) $ equipped with any topology $\tau_2$ in which the maps $e^n_x : f \rightarrow f^{(n)} (x)$ are continuous for every $n \in \N$, $x \in [0,1]$. The finite-open topology $\tau_{fo}$ on the dual spaces $V^*_i$ is given by the subsets $(\{ p \}, O) \subset V_i^*= CL(V_i, \R)$ as a subbasis.  

\begin{Theorem}
\label{sehrbloed}
For $i=1,2$, every metrizable topology on $V_i^* $ is strictly coarser than the finite-open topology.
\end{Theorem}

\v {\bf Proof.} The proof consists of two parts: in the first one, we invoke a theorem from \cite{hhK} to show that for every continuous map $A: V_i \rightarrow (V_i^* , \tau_{fo})$, the map $\tilde{A}: V_i \times V_i \rightarrow \R$ given by $\tilde{A}(v,w) := A(v) (w)$ is continuous. In the second one, for every metrizable topology $\tau_m$ on $V_i^*$ we construct a map $A: V_i \rightarrow (V_i^*, \tau_M) $ for which $\tilde{A}$ is not continuous. Now if $\tau_M$ were finer than $\tau_{fo}$ then $A: V_i \rightarrow (V_i^*, \tau_{fo})$ would be continuous as well and the theorem from the first part would apply, proving the continuity of $\tilde{A}$, a contradiction.

\bigskip

\v {\bf First part:} The theorem from \cite{hhK} (there Lemma 0.1.4.) reads:

\begin{Lemma} 
Let $E$ be a metrizable and barrelled l.c.s., $F$ an arbitrary l.c.s. and $n \in \N$. If $X$ is a metrizable topological space and if $g: X \rightarrow L^n_s (E,F)$ is a continuous function, then the map $\tilde{g} : X \times E^n \rightarrow F   $, associated to $g$, is continuous.
\end{Lemma}

\v We apply this to $E= V_i = X$ which is Fr\'echet and therefore barrelled (and locally convex and metrizable), and $F= \R$. We put $n=1$, then we have that for every $g: V_i \rightarrow L^1_s(V_i,W) $ continuous, $\tilde{g}: V_i \times V_i \rightarrow \R $ is continuous again. By definition on p. 14 we have $L^1 (E,F) = L(E,F)$, and the topology is defined on p.14/15 as the one of simple convergence, that is, equicontinuous convergence on finite sets, which corresponds to the finite-open topology. The definition of $\tilde{g}$ appears on p. 17 (Lemma 0.1.2.).

\bigskip

\v {\bf Second part:} Let $\tau_m$ be any metrizable topology on $V^*$ and let $D$ be a metric compatible to $\tau_m$. Then 

\begin{Lemma}
For all $\d, \e >0$, there is $f^i_{\d, \e} \in V_i^*$ with $f_{\d, \e} \in B^{V_i^*}_{\d} (0)$, but $f^i_{\d, \e} (B^{V_i}_{\e} (0)) = \R$.
\end{Lemma}

\v {\bf Proof of the lemma.} In the case $i=1$ consider the continuous linear maps $d_n$. Suppose there is a ball $B_{\e} (0)$ in which all of them are bounded. Then let $M_n := max \{ \vert d_n (B_{\e} (0)) \vert , n \}$ and consider $v \in V_1$ defined by $v_n := M_n^2$. Then there is no $t >0$ with $t v \in B_{\e} (0)$, a contradiction to the assumption that $\tau_1$ is compatible with scalar multiplication. In the case $i=2$ proceed analogously but replace $d_n$ by $\delta_n: f \mapsto f^{(n)} (1- 2^{-n})$ and define $v$ as the locally finite sum of smooth functions $v_n$ with $v_n ^{(n)} (1- 2^{-n}) = M_n$ and $supp(v_n) \subset [1- 2^{-n} - 2^{-n-2 } , 1- 2^{-n} + 2^{-n-2}]$. Then the supports are disjoint and the sum is defined as a smooth function, and there is no $t>0$ with $ t v \in B_{\e} (0)  $. Now if we require additionally that $\vert \vert v_n \vert \vert_{C^k ([0,1], \R)} < 2^{-n}$ for all $k <n$ then $v$ extends even to $1$ in every $C^k$ (and therefore in the smooth) sense. 

\v {\bf Remark:} We could instead of requiring {\em all} the $e^n_x$ to be continuous only require the $e^n_{1/2} $ to be continuous and then translate the $v_n$ as above to $1/2$. Their sum converges with the same arguments. In general, the arguments given above imply that given a sequence of points $p_n$ in $[0,1]$ and a sequence of numbers $a_n$, one can find a smooth function $f$ on $[0,1]$ with $ f ^{(n)} (p_n) = a_n$. 


\bigskip

\v Now for natural $n$ put $f_n := f_{2^{-n-1}, 2^{-n}}$, then we have $f_n \in B_{2^{-n-1}}^{V_i^*} (0)$ and $f_n (B_{2^{-n}}^{V_i} (0)) = \R$, thus we can pick $w_n' \in B_{2^{-n}}^{V_i} (0)$ with $f_n (w_n') > n $, and there are real numbers $s_n \geq 1$ with $w_n := s_n w_n ' \in B_{2^{-n}}^{V_i} (0) \setminus B_{2^{-n-1}}^{V_i} (0)$ and still $f_n (w_n) >n$. As $V_i$ is a metric space and therefore paracompact, for each $n$ there is a continuous function $\psi_n \in C^0 (V_i, \R) $ with $\psi_n (w_n ) = 1$ and $supp (w_n ) \subset B_{2^{-n-2}} (w_n) \subset V_i \setminus B_{2^{-n-2}} (0)$. The triangle inequality implies that $supp (\psi_i) \cap supp(\psi_{i+4}) = \emptyset $, thus we can define $f:= \sum_{i \in 4 \N} \psi_n \cdot f_n$, and this is a continuous function as every point in $V_i$ has a neighborhood in which every but one term of the sum vanishes. Thus $f \in C^0 (V_i, (V_i^*, \tau_m))$ ($f$ is even a sub-isometry if $D$ has starshaped balls, but we will not need this fact). But $\tilde{f}: V_i \times V_i \rightarrow \R$ is {\em not} continuous, not even $\hat{f} := f \o \Delta : V_i \rightarrow \R$, where $\Delta: V_i \rightarrow V_i \times V_i$ is the diagonal map, so $\hat{f} (v) := f(v) (v)$. This is because $\hat{f} (v_n) = f(w_n) (w_n) > n $, while $w_n \rightarrow 0$.

\bigskip

\v Now we put together the two parts: If there is a metrizable topology on $V_i^*$ and it is finer than $\tau_{fo}$, then we apply the first part to $f$ and conclude that $\tilde{f}$ is continuous, a contradiction.  \hfill \qed

\section{Restriction to tame linear maps} 

\v Another idea is to not consider all of $CL(V,W)$ but only a part $A$ of it. Our first and only try is the space of tame linear maps $TL(V,W)$ which we want to introduce now. The interest in them stems from the fact that every differential operator of degree $k$ corresponds to a $k$-tame map in the natural metrics on spaces of sections (cf. \cite{rsH}, \cite{oM}). A {\bf pre-Fr\'echet space} is a locally convex metric vector space. Consider two pre-Fr\'echet spaces $V,W$. Let $U \subset V$ be open. A map $f: U \rightarrow G$ is called {\bf tame} if for every $u \in U$ there is a neighborhood $A$ of $u$ and $r,b \in \N$ and  $C_n \in \R$ such that for all $a \in A$ and all $n \geq b$ we have 

$$ \mu^W_n (f(a) - f(u)) \leq C_n  (1 + \mu^V_{n+r} (a-u)) , $$ 

\v where the $\mu$ are the respective Minkowski functionals. Now it is esy to see that tameness implies continuity. Less easy to see is the following theorem which gives more restrictive conditions for {\em linear} tame maps: 

\begin{Theorem}[cf. \cite {rsH}]
Let $V, W$ be pre-Fr\'echet spaces. Then any linear $f: V \rightarrow W$ is tame if and only if there are $r,b \in \N$ and $K_n \in \R$ such that for all $v \in V$ we have

$$ \mu_n (f(v)) \leq  K_n \cdot \mu_{n+r} (V) $$

\v for all $n \geq b$. In this case we call $f \in CL(V,W)$ {\bf $r$-tame with basis $b$}. 
\end{Theorem}

\v {\bf Proof.} One direction is trivial. For the other one, assume $f$ is tame. Then for some neighborhood $U$ of $0$ we have

$$ \mu_n (f(v)) \leq C_n (1 + \mu_{n+r} (v)) $$

\v for all $n \geq b$, $v \in U$. Now look for $B \geq b$, $\e >0$ with $\{ v \vert \mu_{B+r} (v) \leq \e \} \subset U$. Choose $v \in V \setminus \{ 0 \} $ and put $g:= \e v / \mu_{B+r} (v) $. Then $\mu_{B+r} (g) = \e$ and 

\bea
\label{nonlintame}
\mu_n (f(g)) \leq C_n (1 + \mu_{n+r} (g)).
\eea 

\v Linearity of $f$ implies $f(g) = \e f(v) / \mu_{B+r} (v)$. Plugging this in into \ref{nonlintame} yields

$$ \frac{\e}{\mu_{B+r} (v)} \cdot \mu_n (f(v)) \leq C_n (1 + \frac{\e}{\mu_{B+r} (v)} \mu_{n+r} (v))  ,$$

\v so for $n \geq B$ we get

\bean
\mu_n (f(v)) &\leq& C_n (\frac{\mu_{B+r} (v)}{\e} + \mu_{n+r} (v))\\
&\leq& C_n (\frac{\mu_{n+r} (v)}{\e} + \mu_{n+r} (v) ) = C_n (1 + \frac{1}{\e}) \mu_{n+r} (v) =: K_n \mu_{n+r} (v), 
\eean

\v which shows that $f$ is $r$-tame with basis $B$. \hfill \qed

\bigskip

\v Obviously, $f \in CL(V,W)$ is $r$-tame with basis $b$ if and only if $f(B(M+r)) \subset C_M \cdot B(M)  \qquad \forall M > b$ for $B(N) := B_{2^{-N}}$ and the $C_N$ being arbitrary real constants. These maps form a subspace $T_{r,b} L(V,W) \subset CL(V,W)$. If there is a natural $b$ such that $f$ is $r$-tame with basis $b$ we call $f$ {\bf $r$-tame} and collect these maps to the space $T_rL(V,W)$. The composition of two tame maps is easily seen to be tame again, but the order of tameness adds up: If $f \in T_rL(V,W)$ and $g \in T_sL(W,X)$, then $f \o g \in T_{r+s} L(V,X)$. Therefore, if we are interested in forming algebras of linear maps or for some other purpose, it seems desirable to collect all tame linear maps in one space, irrespective of their tameness order. Thus we define the space of tame maps $TL(V,W) = \bigcup_{r \in \N} T_r L(V,W) $. Obviously $B_rL(V,W)$ as defined in \cite{oM} is contained in $TL(V,W)$, but not in a single $T_rL(V,W)$. Moreover, we have:

\begin{Theorem}[cf. \cite{rsH}]
\label{tamenessiscontinuity}
Let $V,W$ be pre-Fr\'echet spaces and let the metric of $V$ or of $W$ be tamely equivalent to a norm, then, for all $r,b \in \N$, we have $CL(V,W) = T_{r,b} L(V,W)$, and in particular, $TL(V,W) = CL(V,W)$. Therefore for linear maps between normed spaces tameness is continuity.    
\end{Theorem}

\v {\bf Proof.} Let $f: V \rightarrow W$ be linear and continuous. If $V$ is normable with a norm $\nu$, then for all $n \in \N$, there is an $\e >0$ with $f(\e \cdot B_1^{\nu} (0)) = f(B_{\e}^{\nu} (0)) \subset C(n)$. If $W$ is normable with a norm $\nu$, then there is an $m \in \N$ with $f(c(m)) \subset B_1^{\nu} (0) = 2^m  C(m)$. \hfill \qed

\bigskip

\v {\bf Remark:} Also if both $V$ and $W$ are normed spaces, not every {\em nonlinear} continuous map is tame, consider the continuous function $x \mapsto x^{1/3} $ on $\R$. But at least $C^1 $ maps between finite-dimensional normed spaces are easily seen to be tame.

\bigskip

\v One can now try to replace continuity by tameness in the foundational theorems of infinite-dimensional analysis. For example, it is easy to see that 

\begin{Theorem}
In every pre-Fr\'echet space its metric is a tame (nonlinear) function.  \hfill \qed
\end{Theorem}

\v and, by composing the metric with an appropriate function $f \in C^1([0, \infty), [0,1])$ to conclude that

\begin{Theorem}
\label{partitions}
Every pre-Fr\'echet space has tame partitions of unity. \hfill \qed
\end{Theorem}

\v The subspaces $T_{r,b}L(V,W)$ can be metrized by a quite natural Fr\'echet metric which corresponds to the minimal choice $c_M$ of the $C_M$ if all balls in $W$ are compact and which consists in the familiar Fr\'echet metric for the real sequence $\{ \vert \vert A \vert \vert_M:= \mu_{M} (A (M+r)) \}_{M \geq b}$. 

\begin{Theorem}
The space $T_{r,b} L(V,W)$ with the metric above is a Fr\'echet space.
\end{Theorem}

\v {\bf Proof.} We have to show completeness only. So let $\{ A_n \}_{n \in \N}$ be a Cauchy sequence in the metric, then it is a Cauchy sequence in every $\vert \vert \cdot \vert \vert_M$. Thus the values $A_n (v)$ for a fixed vector $v$ form a Cauchy sequence in every Minkowski functional of $W$. Therefore completeness of $W$ implies that they converge to a point $A(v)$. Then the map $A$ defined as pointwise limit is linear and continuous by the usual arguments, and every $\vert \vert  A \vert \vert_M $ is finite, again because the $\vert \vert A_n \vert \vert_M$ form a Cauchy sequence. \hfill \qed   

\bigskip

\v Now, as the inclusions $T_{r,r}L(V,W) \subset T_{r+1, r+1} L(V,W)$ are strict for all interesting cases (e.g. $V,W$ spaces of sections of fiber bundles), a result by Narayanaswami and Saxon (\cite{NS}) about direct limits of metric vector spaces shows that there is no way to define a Fr\'echet topology on $TL(V,W)$ with all $T_{r,r}L(V,W)$ closed in $TL(V,W)$. What {\em can} be done, however, is, on the one hand, define a (non-complete) locally convex metric tvs structure on $TL(V,W)$ with all $T_{r,r} L(V,W)$ closed, or, on the other hand, define a Fr\'echet structure on $TL(V,W)$ with at least some $T_{r,r}L(V,W)$ not closed. In the light of the inverse function theorems of Nash and Moser the second way seems to be by far more desirable, so we will follow this approach. To this aim, let us introduce some more non-standard terminology: We call a metric vector space $(V,d)$ {\bf strict} iff for all $v \in V$ we have that $S(v) := \sup_{r >0} d(r v, 0)/r  <  \infty  $. An elementary calculation shows that $S(v)= \ov{\lim}_{r \rightarrow 0 } d(rv, 0) /r$ and that $S(\l v ) = \l S(v)$ for $\l >0$. A counterexample to strictness is provided by $(\R, d) $ with $d(x,y) := \sqrt{\vert x-y \vert}$. We call the metric vector space {\bf s-differentiable} iff the function $m_v: s \mapsto d(sv, 0)$ satisfies $m_v \in C^1([0, \infty))$. Obviously, any differentiable metric vector space is strict. The pull-back of an s-differentiable resp. strict metric by a map which is differentiable along rays is s- differentiable resp. strict (the main example is the map $\g \rightarrow \{  \vert \vert \g \vert \vert_n  : n \in \N \}$ for some seminorms $\vert \vert \cdot \vert \vert_n$ as the seminorms are homogeneous, thus differentiable along rays). We will later see that, unfortunately, most common Fr\'echet spaces are {\em not} strict.


\bigskip

\v Let $V,W$ be pre-Fr\'echet spaces. Local convexity implies that the $c(n) := \ov{conv} (B^V_{2^{-n}} (0))$ and $C(n) := \ov{conv} (B^W_{2^{-n}} (0))$ form zero neighborhood bases. Now for a subset $S$ of $W$ put 

$$\mu_n (S) := \inf \{ r \in \R^+ \vert S \subset r \cdot  \ov{conv}(B^W_{2^{-n}}) (0) \}  = ( \sup \{ s \in \R^+ \vert s \cdot S \subset \cdot \ov{conv}(B^W_{2^{-n}}) (0) \} )^{-1}. $$

\v Now we define $\vert \vert A \vert \vert_{m,n} := \mu_n (A (c(m))) $. Theorem \ref{scalbound} tells us that the scalar multiplication with $N$ for $N \in \N $ is bounded by $N$, therefore applying this for $N=2$ we get $2 \cdot B_r(0) \subset B_{2r} (0)$ and $2 \ov{conv} (B_r(0)) \subset \ov{conv} (B_{2r} (0))$ in both $V$ and $W$. Thus $\mu_{i+1} (S) \geq 2 \mu_i (S)$ for every subset $S$ and therefore

\begin{Lemma}
\label{gradus}
$\vert \vert A \vert \vert_{m, n+1} \geq 2 \vert \vert A \vert \vert_{m,n}  $ and $\vert \vert A \vert \vert_{i+1, j} \leq \frac{1}{2} \vert \vert A \vert \vert_{i,j} $. \hfill \qed
\end{Lemma}

\v As a corollary, we get 

\begin{Proposition}
If a continuous linear map $A: V \rightarrow W$ is $r$-tame with basis $b$, it is $(r+b)$-tame with basis $0$. \hfill \qed
\end{Proposition}

\v Now, for some $a_{i,j} \in \R$, we define 

$$K^a_{i,j} := \{ A \in TL(V,W)  : \vert \vert A \vert \vert_{i,j} < a_{i,j} \} = (c(i) , a_{ij} C(j) ) = (a_{ij}^{-1} c(i) , C(j)), \qquad K^a_j := \bigcup_{i=1}^{\infty} K_{i,j} .$$ 

\v Now for $a_{m+1,n} \geq  a_{m,n}/2$, $K^a_j$ is an {\em ascending} union of convex sets because of the Lemma \ref{gradus}, thus it is convex, in particular circled. Now we choose, for $a_{m,n} := m ^{-n}$, $K_{m,l} := K^a_l$ and $K_l := K_{2,l}$. Then the lemma implies $K_{n+1} \subset K_n$. For a real $m$, the topology generated by $ \{ K_{m,l} \vert l \in \N \}$ we denote by $\tau_l$ and put $t:= \tau_2$. Thus $K_{i,j} := (2^i c(i), C(j)) $.

\begin{Theorem}
Let $V,W$ be metric vector spaces, let $V$ be strict. The $K_j$ and their geometric multiples $2^{-n} K_j$ form the countable base of a Hausdorff tvs topology $t$ on $TL(V,W)$ (which is therefore metrizable) that is coarser than the above topology on any $T_{r,b}L(V,W)$, i.e., if $A_n \rightarrow A$ in $T_{r,b}L(V,W)$, then $A_n \rightarrow A$ in $TL(V,W)$.
\end{Theorem}

\v {\bf Proof.} We have to show that 

\begin{enumerate}
\item{The $K_j$ are point-absorbing,}
\item{For every $n \in \N$ there is an $m \in \N$ with $K_m + K_m \subset K_n$.}
\item{For every $v \in V \setminus \{ 0 \}$ there are $j, n \in \N$ with $n \cdot v \notin  K_j$.}
\end{enumerate}

\v Now let $A \in TL(V,W)$, then there is an $r \in \N$ with $A \in T_{r,r}L (V,W)$. Therefore, for a given $j$, look for $i \in \N$ with $i-j  \geq 2r$, then $K_j$ contains $\{ B \in TL(V,W) : \vert \vert B \vert \vert_{i,j} < 2^{-i } \}$. But $\vert \vert A \vert \vert_{i,j} < \infty$, and the norm is homogeneous, so by scaling, $ \l \cdot A \in K_j$. As for the second feature, we can show that $m= n+1$ works: If $v, w \in K_{n+1}$, then there are $i,j \in \N$ with $\vert \vert v \vert \vert_{i,n+1} < 2^{-i}$ and  $\vert \vert w \vert \vert_{j,n+1} < 2^{-j}$. W.r.o.g. let $j \geq i$, then with the second part of Lemma \ref{gradus} we get $\vert \vert v \vert \vert_{j,n+1} < 2^{-i} 2^{-j+i} = 2^{-j}$ and therefore with the first part of Lemma \ref{gradus} we get $2 \vert \vert v+w \vert \vert_{j,n} < \vert \vert v+w \vert \vert_{j,n+1} < 2 \cdot 2^{-j}  $, and the claim follows. For the last property in the list above let $w \in A(V) \setminus \{ 0 \}$. Then there is an $I \in \N$ with $w \notin C(I)$. We choose a $j \geq I$ and want to show that there is an $n \in \N$ with 

$$nA \notin K_j.$$ 

\v Choose a $v \in A^{-1} (w)$. Now, as $V$ is strict, there is a $S(v) >0$ with $d(\l v, 0) < S(v) \cdot \lambda $ for every $\l >0$. For any natural $N>S(v)$ we have $d(\l v, 0) < N \cdot \lambda$ as well for every $\l >0$ and therefore $d(s N^{-1} v, 0 ) < s$ for any $s= \l N >0 $. In particular, for any natural $i$, we have $d(2^{-i} N^{-1} v, 0) < 2^{-i}$ or, in other words, $2^{-i} N^{-1} v \in B^V_{2^{-i}} (0) \subset c(i)$. As $w= A(v)$, we have $2^i N A(c(i)) \not\subset C(j)$ for any natural $i$, so $\vert \vert NA \vert \vert_{i,j} \geq 2^{-i}$ for all natural $i$, thus $N \cdot A \notin K_{i,j}$ for any natural $i$, which means that $N \cdot A \notin K_j$. Thus the topology is Hausdorff, the basis is countable, therefore the topology is metrizable. For the last statement of the theorem note that if a sequence $a_n$ converges in $T_{r,b}L(V,W)$ to $a$, then for every $j$, $a_n -a$ lies finally in $K_{r+j, j}$, so the sequence converges in $TL(V,W)$. \hfill \qed       

\bigskip


\begin{Theorem}
\label{bestpossible}
Let $V,W$ be metric vector spaces, $V$ strict. Let $a: \N^2 \rightarrow \R$ be a map such that the system of $K^a_l$ defined as above is a a zero neighborhood basis of a vector space topology $\tau^a$ on $TL(V,W)$. Then $\tau$ is coarser than $\tau_3$.
\end{Theorem}

\v {\bf Proof.} As vector addition is continuous, for all $j \in \N$ there is a $J \in \N$ with $K_J + K_J \subset K_j$. Now by deleting some $K_l$'s from the basis and renumbering we get a zero neighborhood basis for the same topology $\tau$ with the property that $K_{j+1} + K_{j+1} \subset K_j$ for all natural $j$. This implies that for all natural $i, k$ there is a natural $I(i,k)$ with $K_{i,j} + K_{k,j} \subset 2 K_{I(i,k), j}$. Then we have $conv (K_{i,j}, K_{k,j}) \subset K_{I(i,k), j}$. Now by setting $i(1) := 1$ and $i(n+1) := I(n, i(n))$ we get $K_j = \bigcup_{n \in \N} K_{i(n), j } $ where the $K_{i(n), j} =: \tilde{K}_{nj}$ now form an {\em ascending} union of convex sets. In particular we have 

$$\{ \vert \vert \cdot \vert \vert _{i(n) , j} < a_{i(n), j} \}  = K_{i(n), j} \subset K_{i(n+1), j  } =   \{ \vert \vert \cdot \vert \vert _{i(n+1) , j} < a_{i(n+1), j} \}  .$$

\v The lemma tells us that the left-hand side is contained in $\{ \vert \vert \cdot \vert \vert_{i(n+1), j} < 2^{- D(m) } a_{i(n) , j} \}$ for $D(m) := i(m+1) - i(m)$, so it is clear that $a_{i(n+1), j} > 2^{-D(n)} a_{i(n), j} > 3^{-D(n)} a_{i(n), j}$ is a sufficient condition for this. But using $\mu_k (v) \leq D(v) \cdot 2^k$ it is easy to see that for all $i \in \N$ there is a $I \in \N$ with $\mu_i (v) > 3^{- (I-i)} \mu_I (V)$, thus it is also necessary, so by filling out the gaps between the $i(k)$ we get $K_j \supset \bigcup_{i \in \N} \{ \vert \vert \cdot \vert \vert_{i,j} < a_{1j} \cdot 3^{-i}   \} $. The rest is scaling. \hfill \qed

\bigskip

\v Right from the definition of $TL(V,W)$ as the union of the $T_{r,r}L(V,W)$ it is quite clear that if all inclusions $T_{r,r} L(V,W) \subset T_{r+1, r+1} L(V,W)$ are proper, the former space cannot be complete (choose a diagonal Cauchy sequence). Now we consider the completion $\ov{TL} (V,W)$ of $TL(V,W)$. An element in the completion we call {\bf almost tame}:  

\begin{Theorem}
Let $V,W$ be metric Fr\'echet spaces, $V$ strict. Then $\ov{TL} (V,W) \subset CL(V,W)$ is a metric Fr\'echet space again.
\end{Theorem}


\begin{Theorem}
If $V,W$ are Banach spaces, then $\ov{TL} (V,W) = CL(V,W)$ with the topology coming from the usual operator norm which corresponds to the topology generated by the palette of s-bounded sets.
\end{Theorem}

\v {\bf Proof.} This is easy to see, as in this case the inequalities of Lemma \ref{gradus} are equalities and therefore the sets $K_j$ are balls in the operator norm. However, note that the {\em metric} on $\ov{TL} (V,W)$ defined this way does {\em not} come from a norm (because of the nonlinear $\Phi$). \hfill \qed

\bigskip

\v Now we define a subset $S$ of a metric vector space to be {\bf $\a$-tame} if there is a $D>0 $ with $\mu_n (S) < D \cdot \a^n$ for all natural $n$, it is called {\bf tame} if it is $2$-tame. It is easy to see that if a subset $A$ is tame (so for some $D>0$, $A \subset D 2^{n } C(n)$) then $S$ as in the definition of strictness is bounded by $2D$ on $A$. The image of a tame curve is in general not tame: Suppose that $S$ as in the definition of strictness is unbounded in every ball of the space $V$ (a property which is easy to check for all spaces of sections with standard metrics), choose $v_n \in B_{2^{-n} (0)}$ with $S(v_n) > 2^{-n}$ and join the $v_n$ by straight line segments on the intervals $[2^{-n-1}, 2^{-n}]$ and extend continuously to $c(0) := 0$. The so defined curve $c$ is a subisometry and therefore tame, but obviously its image is not a tame subset of $V$ as $S$ is not bounded on it. 

\v The family of tame subsets forms a palette $\mathcal{T}$ with $\mathcal{F} \subset \mathcal{T} \subset \mathcal{C}$, and the associated palette topology on $CL(V,W)$ we call {\bf tame-open topology}. We say a pre-Fr\'echet space $F$ to satisfy the {\bf Arzela-Ascoli property} iff for any real sequence $a_n$, the set  $\bigcap_{i=1}^{\infty} \{ f \in F \vert \mu_i (f) < a_i \}  $ is compact. The usual spaces of smooth sections do have this property because of the Arzela-Ascoli theorem. This property is a genuine property of metric vector spaces in the sense that, obviously, a normable space has the Arzela-Ascoli property if and only if it is finite-dimensional. We call a pre-Frechet space {\bf $\a$-full} if it contains a compact non-$\a$-tame subset. It is obvious that for pre-Fr\'echet spaces with the Arzela-Ascoli property, every tame set is compact, therefore the tame-open topology is coarser than the compact-open topology. The following theorem shows that it is strictly coarser in case that $V$ has the Arzela-Ascoli property and is full:



\begin{Theorem}
\label{naja}
The metrizable tvs topology $t := \tau_2$ of $\ov{TL}(V,W)$ is finer than the tame-open topology, but coarser than the bounded-open topology.
\end{Theorem}

\v {\bf Proof.} As we need a lemma from general topology in a slightly more general form than the usual one, let us recall it shortly:

\begin{Lemma}
(i) Let $X$ be a set and $\tau, \tilde{\tau}$ two topologies on $X$, then $\tau \subset \tilde{\tau}$ if and only if, for all $x \in X$, for all $A \in N^{\tau} (x)$ there is a $B \in N^{\tilde{\tau}} (x)$ with $B \subset A$.

(ii) If $X$ has the structure of an abelian group and the topologies are compatible with the group structure, then $\tau \subset \tilde{\tau}$ if and only if, for {\em some} $x \in X$, for all $A \in N^{\tau} (x)$ there is a $B \in N^{\tilde{\tau}} (x)$ with $B \subset A$. \hfill \qed
\end{Lemma}

\v Now, by Lemma \ref{gradus}, we have $2^i C(i) \supset 2^{i+1} C(i+1)$ and therefore, for every subset $A$ of $W$, $(2^i C(i), A) \subset (2^{i+1} C(i+1), A)$. Now we use the set-theoretic fact 

$$(L_I, N) \subset \bigcup_{i \in I} (L_i, N) \subset (\bigcap_{i \in I} L_i, N) $$ 

\v to write 

$$(c(I), 2^{-I} s \cdot C(j)) \subset s \cdot K_j = s \cdot lim_{i \rightarrow \infty} (2^i \cdot c(i), C(j))  \subset  (\bigcap_{i=1}^{\infty} 2^i \cdot c(i), s \cdot C(j)),$$

\v and the claim follows by staring at this line: in the center we have a general element of the zero neighborhood basis of $t$. On the left-hand side there is an element of $\tau_b$. Finally, for every element $Y$ of the zero-neighborhood basis of $\tau_t$ we can find a subset $U$ of the form on the right-hand side with $U \subset Y$. Therefore the bounded-open topology (which is {\em not} a vector space topology but still a vector group topology) is finer than $t$ which in turn is finer than the tame-open topology. \hfill \qed

\bigskip

\v This gives us some tools to handle $t$: For example, if a sequence is Cauchy in the bounded-open topology (which, to stress it again, is {\em not} a tvs topology), then it is Cauchy in $t$ (and therefore converges). And on the other hand, if a sequence converges in $t$, it converges uniformly on tame subsets of $V$. We get another immediate corollary:

\begin{Theorem}
\label{notstrict}
$V_1, V_2$ as in the previous section are not strict.
\end{Theorem}

\v {\bf Proof.} This follows from the theorem of the previous section and of the fact that in strict vector spaces, points are tame, thus the tame-open topology is finer than the finite-open topology. \hfill \qed

\begin{Theorem}
\label{schade}
If $V$ is $\a$-full, $\tau_{\a}$ is not finer than the compact-open topology.
\end{Theorem}

\v {\bf Proof.} We choose a non-tame compact set $K$ which exists because of fullness of $V$ and any open set $O \subset W$. We want to show that $(K,O)$ does not contain any $d$-open set. For this, we have to show that for every natural $j$ there is a map $f_j \in K_j$ but $f \notin (K,O)$. We will choose $f_j := \a_j \cdot w$ for a $w \in W \setminus O$ chosen arbitrarily, but fixed for all $j$. The requirement $f \notin (K,O)$ is then satisfied if for any $j$ we find a $v_j \in K$ with $\a_j (v_j ) = 1$. In the same time, we have to show that $f_j \in K_j = \bigcup_{i \in \N} K_{ij}$, that is, we have to find a natural $i$ with $f(C(i)) \subset 2^{-i} C(j) $, or, equivalently, $\vert \vert \a_j \vert \vert_i < 2^{-i} m_j^{-1} $, where $m := \vert \vert w \vert \vert_j$. In the light of the tame Hahn-Banach theorem (applied to $\R w$ and the sublinear functional $\vert \vert \cdot \vert \vert_i$ and keeping in mind that the normability of $\R$ implies that $CL(V,\R) = TL(V,\R)$) this is the same as showing that $\vert \vert v_j \vert \vert_i < 2^{-i} m^{-1}$. So, in summary, we have to find, for every $j \in \N$, a $v_j \in K$ with $\vert \vert v_j \vert \vert_i > 2^{-i} m_j^{-1}  $ (and then define $\a_j$ correspondingly). But the existence of such a $v_j$ is guaranteed precisely by the $\a$-fullness of $V$.     \hfill \qed 


\bigskip

\v So, put together with Theorem \ref{bestpossible}, we get that if $V$ is $3$-full then  no tvs topology on $TL(V,W)$ of the form $\tau^{\a}$ can be finer than the compact-open topology. This is in contrast to the topologies on the spaces $T_{r,b} (V,W)$:

\begin{Theorem}
\label{good}
The Fr\'echet topology on $T_{r,b}L(V,W)$ is finer than the compact-open topology $\tau_{CO}$, and the evaluation map $eva: T_rL(V,W) \times V \rightarrow W$ is continuous.

\end{Theorem}

\v {\bf Proof.} We consider the usual series of seminorms $\vert \vert A \vert \vert_{i} = \vert \vert A \vert \vert_{j+r, j}$ for $j \geq b$. As we can calculate up to tame equivalence, let w.r.o.g. be the metrics on $V,W$ be of sum form. Let a compact $K \subset V$ and an open set $U \subset W$ be given. Then let $a_n$ be the maximum of $\vert \vert \cdot \vert \vert_n$ on $K$. Let $r >0$ with $B_r^W (0) \subset U$. Choose a natural $i$ with $2^{-i} < r/2$, then

$$ \{ w \in W : \vert \vert w \vert \vert_j \leq \Phi^{-1} (r) \qquad \forall j \leq i \} \subset B_r^W (0) .$$

\v Therefore it is sufficient to show $\vert \vert A(K) \vert \vert_j \leq \Phi^{-1} (r)  $ for all $j \leq i$. This is the case if $\vert \vert A \vert \vert_j a_{j+r} \leq \Phi^{-1} (r)  $, or, equivalently and with $M:= max \{1, max_{k= 1,..., i+r} a_k \}$, $\vert \vert A \vert \vert_j \leq \Phi ^{-1} (r) \cdot M^{-1} $ for all $j \leq i$. This is satisfied if 

$$\sum 2^{-j} \Phi (\vert \vert A \vert \vert_j ) \leq 2^{-i} \Phi (\Phi^{-1} (r) \cdot M^{-1}) =: \e,$$

\v so $B_{\e} (0)$ is contained in $(K,U)$. For the second assertion, let $n \in \N$ be given, then, whenever $2^{-n} \Phi (1) > \d$, we get 

$$ B_{2^{-n}}^W (0)  \supset f(B_{\d}^{T_rL(V,W)} (0) \times B_{2^{-n-r} }^V (0) ).  \qquad \qed $$ 

\bigskip

\v Now, in the light of the remark after Theorem \ref{schade}, let us have a look at whether the standard Fr\'echet spaces are $3$-full. To that purpose, we define a pre-Fr\'echet space $F$ to be {\bf step-full} if, for all $s >1$ there are $M(s) \in \R, v_s \in F$ with $ s^i <   \mu_i (v)  < M(s) \cdot (4s)^i $ for all $i \in \N $. Normed spaces are never step-full as their Minkowski functionals grow as a geometric sequence. As the condition on the vector in the definition of step-fullness is preserved by isometries and as there are the functions $sin ((2s)x)$, we get  

\begin{Theorem}
Let $F$ be a pre-Fr\'echet space and $i: F_0 \rightarrow F$ an isometric linear embedding. Then $F$ is step-full. In particular, all spaces of sections of fiber bundles with the standard sum or sup metrics are step-full. \hfill \qed
\end{Theorem}

\begin{Theorem}
\label{step->number}
Let the step-full pre-Fr\'echet space $F$ have the Arzela-Ascoli property. Then it is $s$-full, for each real $s$.   
\end{Theorem}

\v {\bf Proof.} Consider the set $M(s) \cdot \bigcap (4s)^i C(i)$, for an $s \geq 1$. The Arzela-Ascoli property ensures that it is compact. Step-fullness implies that it is not $s$-tame: Suppose that there is a $D >0$ with $M(s) \bigcap (4s)^i C(i) \subset D \bigcap s^i  C(i) $. Then put $j : = min \{ i \in \N \vert 4^i > D/M  \}$ to see that $v_s$ as in the definition of step-fullness is in $\bigcap (4s)^i C(i)$, but not in $D \bigcap s^i  C(i) $, a contradiction. \hfill \qed 

\bigskip


\v Now we get an analogue of Theorem \ref{sehrbloed} for step-full Fr\'echet spaces:

\begin{Theorem} 
Let $V$ be a step-full pre-Fr\'echet space and $W$ be a pre-Fr\'echet space containing a vector $w$ with $S(w) < \infty$. Let $t$ be any topology on $TL(V,W)$ compatible with scalar multiplication. Then there is a tame map $A: U \rightarrow TL(U,W)$ such that $E_A: U \rightarrow W$ given by $E(u) := A(u) (u)$ is not continuous. In particular, the evaluation map $eva: TL(V,W) \times V \rightarrow W$, given by $eva (A, u) = A(u)$ is not continuous, and $V$ is not strict.
\end{Theorem}

\v {\bf Proof.} We choose $r := \frac{1}{2} d(w,0)$. By putting $A := f \cdot w$, by the assumption that $S(w) < \infty$ we reduce the problem to the case $W= \R$: we have to find $f\in TL(V, \R) = CL(V,\R)$ with $f (v_i ) = i$. First, for all $n \in \N$, we construct a tame linear map $\tilde{f}_n \in CL(V, \R)$ which is unbounded on $B_{2^{-n}} (0)$ by the following lemma:


\begin{Lemma}  
Let $V$ be an step-full Fr\'echet space, then for all $\e >0$ there is an $\tilde{f}_{\e} \in TL(V,\R) = CL(V, \R)$ such that $\tilde{f}_{\e} (B_{\e} (0)) $ is unbounded in $\R$. 
\end{Lemma}

\v {\bf Proof of the lemma.} Let $2^{-i-1} < \e  $, then, by step-fullness, consider $v_n $ with $\mu_j (v_n) < 2^{-2n}$ for all $j < i$ and $\mu_j (v_n ) > 1 $ for $j \geq i$. Then, for every $n $, with the Hahn-Banach Theorem we can construct a continuous linear map $\tilde{f}_{\e}^{(n)} $ with $\tilde{f}_{\e}^{(n)} (v_n) := 2^n $ and $\tilde{f}_{\e}^{(n)} < 2^{-n} \mu_{i-1}$ pointwise. Then it is easy to see that the partial sums $\sum_{n= 1}^N \tilde{f}_{\e}^{(n)}$ define pointwise and the map $\tilde{f}_{\e}$ defined as pointwise limit is continuous and linear. But as all $v_n$ are in $B_{\e} (0)$, this set is mapped on an unbounded set in $\R$ by $\tilde{f}_{\e}$.  \hfill \qed 

\bigskip

\v Now we set $\tilde{A}_n = \tilde{f}_n \cdot w$ and proceed as in the proof of Theorem \ref{sehrbloed}: we set $A_n := t_n \cdot \tilde{A}_n  $ for real numbers $t_n$ chosen in a way that $A_n \in B_{2^{-n}}^{TL}(0)$. Still the $A_n$ are unbounded on the $B_{2^{-n}}^V (0) $. Now choose $v_n \in B_{2^{-n}} (0) \setminus B_{2^{-n-1}} (0)$ with $\vert A_n (v_n) \vert > n$, and let $\psi_i $ be a tame function which takes the value $1$ at $v_i$ and whose support is contained in $B_{2^{-i-1}} (v_i)$ (exists because of Theorem \ref{partitions}), then the function $A:= \sum_{i \in 2 \N} \psi_i \cdot A_i :  U \rightarrow TL (V,W)$ is continuous, but $E_A: v \mapsto A(v) (v) $ is not continuous. For the last statement, observe that $E_A = eva \o (A, \1)$, and then proceed as in Theorem \ref{notstrict}.  \hfill \qed

\bigskip

\v For $V,W$ Fr\'echet spaces, $V$ strict, we define iteratively spaces of tame linear maps by $T_r^1 L(V, W) = T_r L(V,W) $ and $T^{n+1}_r L(V,W) := T_rL(V, T^n_r L (V,W))$, $T^1 L(V,W) := TL(V,W)$ and $T^{n+1} L(V,W) := TL(V, T^n L (V,W))$ as well as $\ov{T}^1 L(V,W) := \ov{TL} (V,W)$ and $\ov{T}^{n+1} L (V,W) := \ov{TL} (V, \ov{T}^n L (V,W))  $. Now we can define at least four different types of almost tame $C^k$ maps betwen a strict and a general pre-Fr\'echet space:

\begin{Definition}
Let $U \subset V$ be an open subset of a strict Fr\'echet space and $W$ a Fr\'echet space, let $f \in T(U,W)$. Then we define $f \in T_r^k (U, W)$ iff  $D^l f: U \rightarrow T^l L(V,W) $ exists and is tame for all $l \leq k$. We define $f \in t^k (U,W)$ iff $d^l f: U \times V^l \rightarrow W$ exists and is tame for all $l \leq k$. We write $f \in T^k(U,W)$ iff $D^l f: U \rightarrow T^l L(V,W) $ exists and is tame for all $l \leq k$, and $f \in \ov{T}^k (V,W) $ iff $D^l f : U \rightarrow \ov{T} ^l L (V,W) $ is tame for all $l \leq k$. 
\end{Definition}

\v One example of a nonlinear map which is $T_r$-smooth (and, of course, $t$-smooth) is $\gamma \mapsto g(\nabla_V \gamma, \gamma )$ for $\gamma $ a section of a Riemannian vector bundle with fiber metric $g$ and metric connection $\nabla$ and $V$ a parallel vector field in the base manifold (cf. \cite{oM} where it is shown that this is a bounded-smooth map). Obviously, for $r<s$ and $k \leq \infty $, we have always inclusions $T_r^k \subset T_s^k \subset T^k \subset \overline{T}^k$. 


\v If $V$ and $W$ are Banach spaces, the spaces $t^k (U,W)$ correspond to the spaces of Michal-Bastiani differentiability, the spaces $T_r^k (U,W)$ to Fr\'echet differentiability. Helge Gl\"ockner showed in \cite{hG} that in this case we have $T^k (U,W) \subset t^k (U,W) \subset T^{k-1} (U,W)$. In the case of non-normable spaces we have:


\begin{Theorem}
For all natural $k$, we have the inclusion $T_r^{k+1}(U,W) \subset C^k (U,W)$ with the natural identification of $d^k f$ and $D^k f$. 
\end{Theorem}

\v {\bf Proof.} This is an easy consequence of a standard inductive argument using the continuity of the evaluation map in $T_rL(V,W)$ as established in Theorem \ref{good}. \hfill \qed

\bigskip

\v Resuming, in typical cases, the $T_rL$ spaces between metric vector spaces can be given a useful metrizable topology while the $TL$ spaces cannot. This can be useful in the treatment of spaces of differential operators of bounded degree, but it does not allow for the treatment of {\em algebras} of such operators unless of algebraic ones. Now, in the usual way we can define almost tame manifolds, $T_r^k$ maps between $T_r^k$ manifolds etc. in exactly the same way we define $C^k$ manifolds. In a subsequent analysis we will explore further the relations between the different notions of differentiability, give an exponential map theorem like the one for continuous maps in the compact-open topology, as well as establish two inverse function theorems for these maps.  

\bigskip







\v Finally, I want to give a list of open questions which may have some importance:

\begin{enumerate}
\item{Are there interesting strict vector spaces?}
\item{Is there another definition of $\ov{TL} (V,W)$ not as a completion of something else, but intrinsically as a subspace of $CL(V,W)$?}
\item{Is the metric topology on $TL(V,W)$ coarser than sbounded-open topology if $V$ is not normable?}
\item{Are the metric topologies on $TL(V,W)$ and $T_rL(V,W)$ generated by some palettes? If so, by which ones?}
\item{Are there nice examples of non-normable Fr\'echet spaces which are not $s$-full for some $s>0$?}
\item{Is there a simple criterion to decide when $TL(V,W)$ is strict?}
\item{Is there a sufficiently rich class $C$ of metric vector spaces and continuous linear maps between them such that for the latter ones form again an element of $C$? $TL$ satisfies this statement only almost because of the requirement of strictness on the first space. }
\item{We have seen that no vector space topology on $CL(V, N)$ for $N$ normable makes $eva$ continuous. Is there a compatible metric on $T_rL(V,W)$ that makes $eva$ tame? This would prove $T_r^k (U,W) \subset t^k(U,W)$.}
\end{enumerate}

\begin{small}

\end{small}

\end{document}